\documentclass[12pt]{article}

\usepackage{amssymb}

\author{Alexander L. Zuyev}
\title{Stabilization of the spatial oscillations of an elastic system model}
\begin{document}
\maketitle
\begin{abstract}A system of partial differential equations
describing the spatial oscillations of an Euler-Bernoulli beam
with a tip mass is considered. The linear system considered is
actuated by two independent controls and separated into a pair of
differential equations in a Hilbert space. A feedback control
ensuring strong stability of the equilibrium in the sense of
Lyapunov is proposed. The proof of the main result is based on the
theory of strongly continuous semigroups.
\end{abstract}

%%%%%%%%%%%%%%%%%%%%%%%%%%%%%%%%%%%%%%%%%%%%%%%%%%%%%%%%%%%%%%%%%%%%%%%%%%%%%%%%
\section{Introduction}

Dynamical models of flexible-link robot manipulators are generally described by
a set of coupled ordinary and partial differential equations, that gives rise
to series of mathematical control problems in infinite dimensional
spaces~\cite{Fat,KS,LL,LGM,ZuCDC2003,ZuAutomatica}. However, finite dimensional
approximate models obtained by the assumed modes and finite elements methods
are used more frequently for solving the motion planning and stabilization
problems~\cite{CBG,TPK,ZuCDC2006,ZuMPE}. It should be emphasized that the majority of
publications in this area is concentrated on planar manipulator models with a
free end. To study spatial manipulators with a tip mass, the mathematical model
that describes the motion of a multi-link manipulator under the action of
gravity and controls (torques and forces) was proposed in~\cite{ZuIAMM2005}.

The goal of this paper is to study the stabilization problem of
the control system derived in~\cite{ZuIAMM2005} for the particular
case of a manipulator with one flexible link.

\section{Equations of Motion}

A mechanical system consisting of $n$ Euler-Bernoulli beams and a rigid body as
a load was introduced in~\cite{ZuIAMM2005}. In this paper, we assume that $n=1$
and neglect controlled rotations of the load ($\varphi_J=0$ and $c=0$ in the
notations of~\cite{ZuIAMM2005}). Thus, the beam deflection at time $t$ is
defined by functions $y(x,t)$ and $z(x,t)$ in a rotating Cartesian frame, where
$x\in [0,l]$ is the spatial coordinate, $l$ is the length of the beam. The
above Cartesian frame is obtained from the fixed one by subsequent rotations on
the angle $\varphi_T(t)$ (turning angle) and $\varphi_R(t)$ (raising angle).
The system is controlled by torques $M_T$ and $M_R$ applied at the bottom part
of the beam. For each constant value $\varphi_R^0$, there is the control torque
$M_R=M_R^0$ implementing an equilibrium $\varphi_T(t)=0$,
$\varphi_R(t)=\varphi_R^0$, $y(x,t)=0$, and $z(x,t)=z_0(x)$. The linearized
system of differential equation describing oscillations around the equilibrium
can be written as follows (see~\cite{ZuIAMM2005}):
\begin{equation}
\ddot y(x,t)+ \frac{1}{\rho} \Bigl(c_z y''(x,t)\Bigr)'' =
\psi_T(x) \ddot \varphi_T,\;x\in (0,l), \label{PDEy_lin}
\end{equation}
\begin{equation}
\ddot {\tilde z}(x,t)+ \frac{1}{\rho} \Bigl(c_y {\tilde
z}''(x,t)\Bigr)'' = g  \tilde \varphi_R \sin\varphi_R^0 - x \ddot
{\tilde \varphi}_R, \label{PDEz_lin}
\end{equation}
\begin{equation}
\left. y \right|_{x=0} = \left. \tilde z \right|_{x=0} =
0,\;\left. y' \right|_{x=0} = \left. \tilde z' \right|_{x=0} = 0,
\label{BCyz_zero}
\end{equation}
\begin{equation}
\left.\frac{1}{m} (c_z y'')' - \ddot y +\psi_T(x)\ddot \varphi_T
\right|_{x=l} =0, \label{BCy_lin1}
\end{equation}
\begin{equation}
\left. -c_z y'' - J_3 \ddot y' + J_3 \psi_T'(x)\ddot\varphi_T
\right|_{x=l} =0, \label{BCy_lin2}
\end{equation}
\begin{equation}
\left.\frac{1}{m} (c_y \tilde z'')' + g{\tilde \varphi}_R
\sin\varphi_R^0 - \ddot {\tilde z} - l \ddot {\tilde
\varphi}_R\right|_{x=l} =0, \label{BCz_lin1}
\end{equation}
\begin{equation}
\left. c_y \tilde z'' + J_2 (\ddot{\tilde \varphi}_R+\ddot{\tilde
z}')\right|_{x=l} =0, \label{BCz_lin2}
\end{equation}
%%%%
$$
\left\{I_0 +
(I_1+J_1)\sin^2\varphi_R^0+(I_3+J_3)\cos^2\varphi_R^0+\right.
m_0(R-d\cos\varphi_R^0)^2 + m(R-l\cos\varphi_R^0)^2 +
$$
$$ \left.+
\int_0^l (R-x\cos\varphi_R^0)^2\rho\,dx\right\}\ddot \varphi_T+
\int_0^l (R-x\cos\varphi_R^0)\ddot y \rho\,dx +
$$
\begin{equation}
+ \left.\left\{m R \ddot y -(ml\ddot y  + J_3\ddot y')
\cos\varphi_R^0\right\}\right|_{x=l} = M_T, \label{Lagrange_y_lin}
\end{equation}
$$
\left\{I_2 +J_2 + m_0 d^2 + m l^2 + \int_0^l x^2 \rho\,dx
\right\}\ddot {\tilde \varphi}_R + \int_0^l \ddot{\tilde z} x
\rho\,dx + \left.\left\{ml \ddot{\tilde z} + J_2 \ddot{\tilde
z}'\right\}\right|_{x=l} -
$$
$$
 -g\left\{\int_0^l \tilde z\rho \,
dx + \left.m{\tilde z}\right|_{x=l}\right.\left. +\left(m_0d
+ml+\int_0^l x\rho\, dx\right){\tilde
\varphi}_R\right\}\sin\varphi_R^0 -
$$
\begin{equation}
-g\left\{\int_0^l z_0\rho\,dx+m z_0(l)\right\}{\tilde \varphi}_R
\cos\varphi_R^0= {\tilde M}_R, \label{Lagrange_z_lin}
\end{equation}
where $\tilde z(x,t) = z(x,t)-z_0(x)$, $\tilde\varphi_R(t) =
\varphi_R(t)-\varphi_R^0$, $\tilde M_R = M_R - M_R^0$, and
$$
\psi_T(x)=x \cos\varphi_R^0 - z_0(x) \sin \varphi_R^0 - R.
$$
We use dots to denote derivatives with respect to time $t$, and primes to
denote derivatives with respect to the space variable $x$. The procedure for
computing $z_0(x)$ and $M_R^0$ is given in~\cite{ZuIAMM2005}.

The parameters in~(\ref{PDEy_lin})-(\ref{Lagrange_z_lin}) have the following
physical meaning: $\rho(x)$ is the mass per unit length of the beam, $c_z(x) =
E(x)I_z(x)$, $c_y(x) = E(x)I_y(x)$, $E(x)$ is Young's modulus, $I_z(x)$ and
$I_y(x)$ are moments of inertia of the cross section of the beam with respect
to the axes $z$ and $y$, $m$ is the payload mass, $J_1$, $J_2$, and $J_3$ are
central moments of inertia of the payload, $R$ is the platform radius, $I_0$ is
the moment of inertia of the platform, $I_1$, $I_2$, and $I_3$ are moments of
inertia of the hub, $m_0$ is the hub mass, $d$ is the distance between the
origin of the rotating Cartesian frame and the hub center of mass.

To simplify these equations we substitute expressions~(\ref{PDEy_lin}),
(\ref{PDEz_lin}), (\ref{BCy_lin1})-(\ref{BCz_lin2}) for $\ddot y(x,t)$, $\ddot
{\tilde z}(x,t)$, $\left.\ddot y, \ddot y', \ddot {\tilde z},\ddot {\tilde
z}'\right|_{x=l}$ into~(\ref{Lagrange_y_lin}), (\ref{Lagrange_z_lin}) and
perform integration by parts with regard for the boundary
conditions~(\ref{BCyz_zero}). As a result, equations~(\ref{Lagrange_y_lin}) and
(\ref{Lagrange_z_lin}) take the following form:
\begin{equation}
\ddot \varphi_T = u_T,\; \ddot {\tilde \varphi}_R = u_R,
\label{phidynamics}
\end{equation}
where
$$
u_T = \{I_0 + (I_1+J_1)\sin^2\varphi_R^0+m_0(R-d\cos\varphi_R^0)^2
+
$$
$$
+(I_3\cos\varphi_R^0+J_3 z_0'(l)\sin\varphi_R^0)\cos\varphi_R^0
 + \Bigl(m(l\cos\varphi_R^0-R)z_0(l)+\Bigr.
 $$
\begin{equation}
 +\Bigl.\int_0^l (x\cos\varphi_R^0-R)z_0 \rho
\,dx\Bigr)\sin\varphi_R^0\bigr\}^{-1} \times
 \left\{ M_T -
\left.\left(R (c_z y'')'+c_z y''
\cos\varphi_R^0\right)\right|_{x=0} \right\}, \label{uT}
\end{equation}
$$
u_R =\{I_2+m_0d^2\}^{-1}\times\Bigl\{{\tilde M}_R+c_y
\left.{\tilde z}''\right|_{x=0}+ g\left(\int_0^l \tilde z \rho \,
dx+ m \left.\tilde z\right|_{x=l}+m_0d\right)\sin\varphi_R^0+
$$
\begin{equation}
+g\left(\int_0^l z_0 \rho\,dx + m z_0(l)\right){\tilde
\varphi}_R\cos\varphi_R^0\Bigr\}. \label{uR}
\end{equation}

For each $\tilde \varphi_R(t)$, $y(\cdot, t)$, $\tilde z(\cdot,t)$,
formulae~(\ref{uT}) and (\ref{uR}) establish a one-to-one correspondence
between the torques $(M_T, \tilde M_R)$ and angular accelerations $(u_T,u_R)$.
Thus, we may consider $(u_T,u_R)\in {\mathbb R}^2$ as a new control for the
linear system~(\ref{PDEy_lin})-(\ref{BCz_lin2}), (\ref{phidynamics}).

\vspace{4mm} {\bf 3~~~Main Results}

Consider the following linear space
$$
X  = \Bigl\{{\small \left(\begin{array}{c}\eta(\cdot) \\
\zeta(\cdot)
\\ \phi \\ \omega \\ p \\ q \end{array}\right)}\,:\, \begin{array}{l}\eta\in
H^2(0,l),\;\zeta \in L_2(0,l),\\
\eta(0)=\eta'(0)=0,\\ \phi,\omega,p,q\in \mathbb
R\end{array}\Bigr\},
$$
where $H^k(0,l)$ is the Sobolev space of all functions whose generalized
derivatives of order $j=0,1,...,k$ exist and belong to $L_2(0,l)$. For
$$
\xi_1 ={\small \left(\begin{array}{c}\eta_1 \\
\zeta_1
\\ \phi_1 \\ \omega_1 \\ p_1 \\ q_1 \end{array}\right)}\in X\;\; {\rm and} \;\;
\xi_2 ={\small \left(\begin{array}{c}\eta_2 \\
\zeta_2
\\ \phi_2 \\ \omega_2 \\ p_2 \\ q_2 \end{array}\right)}\in X,
$$
the inner product in $X$ is defined by the formula
$$
\left<\xi_1,\xi_2\right>_X = \int_0^l (\eta_1'' (x)\eta_2''(x) +
\zeta_1(x)\zeta_2(x))\,dx +\phi_1 \phi_2 + \omega_1 \omega_2 + p_1
p_2 + q_1 q_2.
$$
It is easy to check that the norm $\|\xi\|_X =
\sqrt{\left<\xi,\xi\right>_X}$ is equivalent to the standard norm
in $H^2(0,l)\times L_2(0,l)\times {\mathbb R}^4$ (see,
e.g.,\cite[Ch.~3]{Mikhailov}), and hence,
$\left(X,\|\cdot\|_X\right)$ is a Hilbert space.

In order to consider an abstract formulation of the boundary value
problem~(\ref{PDEy_lin})-(\ref{BCz_lin2}), (\ref{phidynamics}),
let us introduce the linear operator $A:D(A)\to X$ and the element
$B\in X$ as follows:
\begin{equation}
A:\;\xi= \left(\begin{array}{c} \eta
\\ \zeta \\ \phi\\ \omega \\ p\\ q
\end{array}\right)\mapsto A \xi= \left(\begin{array}{c}
\zeta
\\
-\frac{1}{\rho} (c \eta'')''+\gamma \phi
\\ \omega \\ 0
\\ \gamma \phi+\frac{1}{m}(c\eta'')'|_{x=l}\\ -
\frac{c}{J}\eta''|_{x=l}
\end{array}\right),\;
B= \left(\begin{array}{c} 0 \\ \psi \\ 0 \\
1 \\ \psi (l) \\ \psi'(l)\end{array}\right),
\label{op_AB}
\end{equation}
where the domain of definition of $A$ is
\begin{equation}
D(A) =\Bigl\{\xi\in X\,:\,\begin{array}{l}\eta\in
H^4(0,l),\zeta\in H^2(0,l),\\
\zeta(0)=\zeta'(0)=0,\\ p=\zeta(l),q=\zeta'(l)\end{array}\Bigr\},
 \label{DA}
\end{equation}
functions $c(x)>0$ and $\psi(x)$ are assumed to be of class
$C^2[0,l]$; $J>0$ and $\gamma$ are constants.

Let $\Bigl(y(x,t),\tilde
z(x,t),\varphi_T(t),\tilde\varphi_R(t)\Bigr)$ be a classical
solution of the boundary-value
problem~(\ref{PDEy_lin})-(\ref{BCz_lin2}), (\ref{phidynamics})
with controls $(u_T(t),u_R(t))$ for $0\le t<\tau$, $\tau\le
+\infty$. Defining
\begin{equation}
\xi_T(t) ={\small \left(\begin{array}{c} y(\cdot,t) \\
\dot y(\cdot,t)
\\ \varphi_T(t) \\ \dot\varphi_T(t) \\ \dot y(l,t) \\ \dot y'(l,t) \end{array}\right)}, \;\;
\xi_R(t) ={\small \left(\begin{array}{c} \tilde z(\cdot,t) \\
\dot {\tilde z}(\cdot,t)
\\ \tilde\varphi_R(t) \\ \dot{\tilde\varphi}_R(t) \\ \dot {\tilde z}(l,t) \\ \dot {\tilde z}'(l,t)
\end{array}\right)},
\label{phaseanalog}
\end{equation}
we see that $\xi_T(t)\in D(A)$ and $\xi_R(t)\in D(A)$ for all
$t\in [0,\tau)$. Consider the pair $(A_T,B_T)$ obtained from
$(A,B)$ by placing $\psi(x)=\psi_T(x)$, $c(x)=c_z(x)$, $J=J_3$
$\gamma=0$ in~(\ref{op_AB}). Similarly, let the pair $(A_R,B_R)$
be obtained from $(A,B)$ by plugging $\psi(x)=-x$, $c(x)=c_y(x)$,
$J=J_2$, $\gamma=g\sin\varphi_R^0$. Then the boundary-value
problem~(\ref{PDEy_lin})-(\ref{BCz_lin2}), (\ref{phidynamics}) is
reduced to the following control system:
\begin{equation}
\dot \xi_T = A_T \xi_T + B_T u_T,
\label{operator_horiz}
\end{equation}
\begin{equation}
\dot \xi_R = A_R \xi_R + B_R u_R,
\label{operator_vert}
\end{equation}
where $(\xi_T,\xi_R)$ is the state and $(u_T,u_R)$ is the control. In the
sequel, we treat this control system as an abstract formulation
of~(\ref{PDEy_lin})-(\ref{BCz_lin2}), (\ref{phidynamics}) with $\xi_T,\xi_R\in
X$ and $u_T,u_R\in {\mathbb R}$. We see
that~(\ref{operator_horiz}),~(\ref{operator_vert}) is separated into two parts,
therefore, the stabilization problem may be solved independently for $\xi_T$
and $\xi_R$. The basic result we shall prove is the following

{\bf Theorem~1.} {\em Consider the abstract Cauchy problem on
$t\ge 0$:
\begin{equation}
\dot \xi(t) = A \xi(t) + B u, \label{cs_AB}
\end{equation}
\begin{equation}
\xi(0)=\xi_0\in X,
\label{Cauchy_IC}
\end{equation}
where $A$, $B$ are given by~(\ref{op_AB}),
$$
u = - \frac{1}{\beta}\left\{k \omega+
\left(\alpha-\gamma\left(\int_0^l\rho\psi\,dx+m\psi(l)\right)\right)\phi+\right.
$$
\begin{equation}
+\int_0^l c\eta''\psi''dx+\left.\Bigl(c
\eta''\psi'-(c\eta'')'\psi\Bigr)\right|_{x=0}
-\left.\gamma\left(\int_0^l\rho\eta\,dx+m\eta(l)\right)\right\},
\label{control_AB}
\end{equation}
$\alpha>0$, $\beta>0$ are large enough constants, and $k>0$ is
arbitrary.

Then the Cauchy problem~(\ref{cs_AB}), (\ref{Cauchy_IC})
with~(\ref{control_AB}) is well-posed on $t\ge 0$ (in the sense of
mild solutions), and the feedback control~(\ref{control_AB})
strongly stabilizes the equilibrium $\xi=0$ of the control
system~(\ref{cs_AB}), i.e., for every $\varepsilon>0$, there
exists $\delta=\delta(\varepsilon)>0$ such that, for every
solution of~(\ref{cs_AB})-(\ref{control_AB}),
$$
\|\xi_0\|_X<\delta \Rightarrow \|\xi(t)\|_X<\varepsilon,\;\forall
t\ge 0.
$$
Moreover, if a semitrajectory $\{\xi(t)\}_{t\ge 0}$
of~(\ref{cs_AB}), (\ref{control_AB}) is precompact in $X$ then the
set of its limit points (as $t\to+\infty$) is an invariant subset
of $Z_0 = \{\xi\in X\,:\, \omega = 0\}$. }

{\em Proof.} Consider a quadratic functional on $X$:
$$
2V(\xi)= \alpha \phi^2 + \beta \omega^2 + \int_0^l
\left\{(\zeta-\psi\omega )^2\rho + {\eta''}^2 c\right\}dx +
$$
\begin{equation}
+m\left\{p-\psi(l)\omega\right\}^2 + J \left\{q-\psi'(l)\omega
\right\}^2 -2\gamma\phi\left\{\int_0^l \eta\rho\,dx + m
\eta(l)\right\}. \label{Lyapunov_horiz}
\end{equation}
Let us compute the time-derivative of $V$ along trajectories
of~(\ref{cs_AB}) when $\xi\in D(A)$:
$$
\dot V(\xi) = \left<\nabla V(\xi),A\xi + B u\right> = $$
$$
=\int_0^l \Bigl(c \zeta'' \eta'' - \zeta \cdot (c
\eta'')''\Bigr)dx + \left.\Bigl(p \,(c \eta'')'-q c
\eta''\Bigr)\right|_{l}+
$$
$$
+\Bigl\{ \alpha\phi+\beta u +\int_0^l \Bigl(\psi \cdot (c\eta'')''
-\rho\gamma(\phi\psi+\eta)\Bigr)\,dx+
$$
\begin{equation}
+\left.\Bigl(c \psi' \eta'' - \psi \cdot (c \eta'')'-m\gamma
(\phi\psi + \eta)\Bigr)\right|_{x=l}\Bigr\}\omega. \label{vdot1}
\end{equation}
By performing integration by parts with regard for
conditions~(\ref{DA}), we get:
$$
\int_0^l\zeta \cdot (c \eta'')''dx = \left.\zeta \cdot (c
\eta'')'\right|_{x=0}^l - \int_0^l \zeta' \cdot (c\eta'')' dx =
$$
$$
=\left.\Bigl(\zeta\cdot (c\eta'')' - \zeta' c
\eta''\Bigr)\right|_{x=l}+\int_0^l \zeta'' c \eta'' dx,
$$
$$
\int_0^l \psi \cdot (c \eta'')''dx = \left.\psi \cdot
(c\eta'')'\right|_{x=0}^l - \int_0^l \psi' (c\eta'')' dx =
$$
$$
=
\left.\Bigl(\psi \cdot (c\eta'')' - \psi' c
\eta''\Bigr)\right|_{x=0}^l+\int_0^l \psi'' c \eta_T'' dx.
$$
Let us substitute these formulae into~(\ref{vdot1}) and use
boundary conditions $p=\zeta(l)$, $q=\zeta'(l)$ from~(\ref{DA}).
As a result, the expression for $\dot V$ takes the following form:
$$
\dot V(\xi) = \Bigl\{ \left(\alpha-\gamma\left(\int_0^l
\rho\psi\,dx+m\psi(l)\right)\right)\phi+\beta u +
$$
\begin{equation}
+ \int_0^l c\psi'' \eta'' dx+\left.\Bigl(c \psi' \eta'' -
\psi\cdot (c \eta'')'\Bigr)\right|_{x=0} -\gamma\left(\int_0^l
\rho\eta\,dx+m\eta(l)\right) \Bigr\}\omega. \label{vdot2}
\end{equation}
If $u$ is defined by~(\ref{control_AB}) then formula~(\ref{vdot2})
yields
\begin{equation}
\dot V(\xi) = - k \omega^2\le 0,\quad (k=const>0). \label{vdot}
\end{equation}

The next step is to prove that $V(\xi)$ satisfies the following
estimates
\begin{equation}
M_1 \|\xi\|_X^2\le 2V(\xi)\le M_2 \|\xi\|_X^2 \label{Vest}
\end{equation}
with some constants $0<M_1\le M_2<+\infty$. On one hand, by
exploiting inequalities $(a+b)^2\le 2a^2 + 2b^2$ and $2ab\le
a^2+b^2$ in~(\ref{Lyapunov_horiz}), we obtain
$$
2V(\xi)\le \alpha \phi^2 + \beta \omega^2 + \int_0^l \Bigl(c
{\eta''}^2+2 (\zeta^2+\psi^2\omega^2)\rho \Bigr)dx +
$$
$$
+2m\Bigl(p^2 + \psi^2(l)\omega^2\Bigr)+
 2 J
\Bigl(q^2 + {\psi'}^2(l)\omega^2\Bigr)+
$$
\begin{equation}
+\gamma^2\phi^2 + 2\left(\int_0^l \eta\rho \, dx \right)^2 +
2\left(m\eta(l)\right)^2.
\label{Vsup}
\end{equation}
Then the Cauchy-Schwartz inequality implies
\begin{equation}
\left(\int_0^l \eta\rho \, dx \right)^2 \le  \int_0^l \eta^2 \,
dx\int_0^l \rho^2 \, dx,
\label{CauchyScwartz1}
\end{equation}
\begin{equation}
\eta^2(l)=\left(\int_0^l \eta'\, dx \right)^2\le \int_0^l dx
\int_0^l {\eta'}^2 dx. \label{CauchyScwartz2}
\end{equation}
The functions $\eta(x)$ and $\eta'(x)$ subject to the boundary
conditions $\eta(0)=\eta'(0)=0$ satisfy Friedrichs' inequalities
of the following form~(cf.~\cite[p.~440]{KS}):
\begin{equation}
\int_0^l {\eta}^2 dx \le \frac{l^2}{2}\int_0^l {\eta'}^2 dx\le
\frac{l^4}{4}\int_0^l {\eta''}^2 dx.
\label{Friedrichs}
\end{equation}
By using inequalities~(\ref{CauchyScwartz1})-(\ref{Friedrichs}) we
conclude that
$$
\left(\int_0^l \eta\rho \, dx \right)^2 + m^2\eta^2(l) \le
\int_0^l \eta^2 \, dx\int_0^l \rho^2 \, dx +
$$
\begin{equation}
+ l m^2 \int_0^l {\eta'}^2 dx \le
\frac{l^3}{2}\left(m^2+\frac{l}{2}\int_0^l \rho^2
dx\right)\int_0^l {\eta''}^2 dx.
\label{eta2est}
\end{equation}
Application of this
inequality in~(\ref{Vsup}) yields an estimate $2V(\xi)\le M_2
\|\xi_T\|_X^2$,
$$
M_2 = \max \Bigl\{\alpha+\gamma^2,2m,2J,2\max_{x\in [0,l]}\rho(x),
$$
$$
\beta+2\int_0^l \psi^2\rho \, dx + 2J {\psi'}^2(l)+2m\psi^2(l),
$$
$$
l^3\left(m^2+\frac{l}{2}\int_0^l \rho^2 dx\right) +\max_{x\in
[0,l]}c(x)\Bigr\}.
$$
On the other hand, we see that the inequality $a^2=(a-b+b)^2\le
2(a-b)^2+2b^2$ implies $(a-b)^2\ge a^2/2 - b^2$. By using the
latter together with $-2ab\ge -\varkappa^2 a^2-b^2/\varkappa^2$
($\varkappa\neq 0$) in~(\ref{Lyapunov_horiz}), we get:
$$
2V(\xi)\ge \alpha \phi^2 + \beta \omega^2 +\int_0^l\Bigl(c
{\eta''}^2 + \frac{\rho}{2} \zeta^2 -  \rho\psi^2
\omega^2\Bigr)dx+
$$
$$
+m\Bigl(\frac{p^2}{2}-\psi^2(l)\omega^2\Bigr)+
J\Bigl(\frac{q^2}{2}-{\psi'}^2(l)\omega^2\Bigr)-
$$
$$
-\varkappa^2\gamma^2 \phi^2 - \frac{1}{\varkappa^2}\left(\int_0^l
\eta \rho \, dx + m \eta(l)\right)^2\ge
$$
$$
\ge \left(\alpha-\varkappa^2\gamma^2\right)\phi^2
+\frac{m}{2}p^2+\frac{J}{2}q^2+\frac{1}{2}\int_0^l \zeta^2 \rho\,
dx+
$$
$$
+\Bigl( \beta-\int_0^l \rho\psi^2\,dx - m\psi^2(l)-J
{\psi'}^2(l)\Bigr)\omega^2+
$$
\begin{equation}
+\left\{\min_{
[0,l]}c-\frac{l^3}{\varkappa^2}\left(m^2+\frac{l}{2}\int_0^l
\rho^2 dx\right) \right\}\int_0^l {\eta''}^2 dx.
\label{Vinf}
\end{equation}
We have also used the inequality~(\ref{eta2est}) here.
From~(\ref{Vinf}) we conclude that $2V(\xi)\ge M_1\|\xi\|_X^2$ and
$$
M_1 =
\min\Bigl\{\alpha-\varkappa^2\gamma^2,\frac{m}{2},\frac{J}{2},\frac{1}{2}\min_{x\in
[0,l]}\rho(x),
$$
$$
 \beta
-\int_0^l \rho\psi^2\,dx -m\psi^2(l)- J{\psi'}^2(l),
$$
$$
\min_{ x\in
[0,l]}c(x)-\frac{l^3}{\varkappa^2}\left(m^2+\frac{l}{2}\int_0^l
\rho^2 dx\right) \Bigr\}>0
$$
provided that
$$
\varkappa^2 > \frac{l^3}{\min_{
x\in[0,l]}c(x)}\left(m^2+\frac{l}{2}\int_0^l \rho^2 dx\right),
$$
$$
\alpha> \varkappa^2\gamma^2,\; \beta> \int_0^l \rho\psi^2\,dx
+m\psi^2(l)+ J {\psi'}^2(l).
$$
For the rest of this paper, we assume that constants $\alpha$,
$\beta$, and $\varkappa$ satisfy the above inequalities.

The estimate~(\ref{Vest}) shows that the two norms $\|\xi \|_X$
and $\|\xi\|_V=\sqrt{V(\xi)}$ are equivalent in $X$. Let us write
the closed-loop system~(\ref{cs_AB}) with the control $u$ defined
by~(\ref{control_AB}) as $ \dot \xi=\tilde A \xi$, where $D(\tilde
A)= D(A)$ is dense in $X$. From inequality~(\ref{vdot}) it follows
that the operator $\tilde A$ is dissipative in $X$ equipped with
the norm $\|\cdot\|_V$. Then the Lumer-Phillips
theorem~\cite[Chap.~1.4]{Pazy} implies that $\tilde A$ is the
infinitesimal generator of a $C_0$ semigroup of contractions,
$\{e^{t\tilde A}\}_{t\ge 0}$, on $X$ (with respect to the norm
$\|\cdot\|_V$). It means that the Cauchy
problem~(\ref{cs_AB})-(\ref{control_AB}) has the unique mild
solution $\xi(t)=e^{t\tilde A}\xi_0$, $t\ge 0$, for every
$\xi_0\in X$, and the above solution is classical if $\xi_0\in
D(A)$. As $\{e^{t\tilde A}\}_{t\ge 0}$ is contractive (under an
equivalent renormalization in $X$), then
$$
\|\xi(t)\|_{V}\le \|\xi_0\|_{V},\quad \forall t \ge 0.
$$
This implies, taking into account the estimate~(\ref{Vest}), that
$$
\|\xi(t)\|_X^2 \le \frac{2V(\xi(t))}{M_1}\le
\frac{2V(\xi_0)}{M_1}\le \frac{M_2}{M_1} \|\xi_0\|_X^2.
$$
The above inequality proves strong stability of the equilibrium
$\xi=0$ in the sense of Lyapunov (we may choose
$\delta(\varepsilon)=\varepsilon\sqrt{M_1/M_2}$ in the definition
of stability).

To conclude the proof we apply LaSalle's invariance
principle~\cite{LaSalle,Shest} (cf.~\cite[Lemma~2]{ZuUMJ2006})
with the functional $V(\xi)$: if a semitrajectory
$\{\xi(t)\}_{t\ge 0}$ is precompact then its $\omega$-limit set,
$\Omega(\xi_0)$, is a non-empty and semi-invariant subset of
$\overline{\{\xi\in D(A)\,:\, \dot V(\xi) =0 \}}= Z_0$. $\square$

{\bf Remark.} By applying the formula~(\ref{control_AB}) to
control systems~(\ref{operator_horiz}), (\ref{operator_vert})
separately and using the representation~(\ref{phaseanalog}), one
can write the feedback control proposed as follows:
$$
u_T = - \frac{1}{\beta}\Bigl\{\alpha\varphi_T + k \dot\varphi_T +
\int_0^l c_z y'' \psi''dx +\left.\Bigl(c_z y''\psi'-(c_z
y'')'\psi\Bigr)\right|_{x=0}\Bigr\},
$$
$$
u_R = - \frac{1}{\beta}\Bigl\{\alpha\tilde\varphi_R + k
\dot{\tilde\varphi}_R - \left. c_y \tilde z''\right|_{x=0}+
$$
$$
+ g\left(\int_0^l  (x\tilde\varphi_R - \tilde z)
\rho\,dx+m(l\tilde\varphi_R-\tilde z|_{x=l})\right)\sin\varphi_R^0
\Bigr\}.
$$
To implement these controls in practice, it is sufficient to
compute $u_T$ and $u_R$ depending on the measurements of
$\varphi_T$, $\tilde \varphi_R$, $\dot \varphi_T$, $\dot {\tilde
\varphi}_R$, $y$, $\tilde z$ at each $t\ge 0$, and then apply
formulae~(\ref{uT}), (\ref{uR}) to find torques $M_T$ and $M_R$.
An advantage of this approach is that no information about the
time-derivatives of $y(x,t)$ and $\tilde z(x,t)$ is needed.

\section{Conclusions}

A feedback control has been derived to stabilize the equilibrium
of a differential equation in a Hilbert space that describes the
motion of a flexible beam with a tip mass. Although the main
result of this paper concerns non-asymptotic stability, further
analysis of the limit behavior of controlled trajectories is
possible by means of the invariance principle. The main difficulty
in this direction is to prove that the semitrajectories are
precompact, which is not an easy task in general (see,
e.g.,~\cite{CA,ZuNASU2007}). We do not study the compactness issue
here, leaving it for future work.

This research is supported by the Ministry of Education and
Science of Ukraine through grant of the President of Ukraine for
Young Scientists.

\end{document}